\documentstyle[11pt,epsf] {article}
\newlength{\mytopmargin}
\newlength{\myleftmargin}
\newtheorem{lemma}{Lemma}
\newtheorem{prop}[lemma]{Proposition}
\begin{document}
\vspace{4cm}
\noindent
{\bf RANDOM  WALKS AND RANDOM FIXED POINT FREE INVOLUTIONS}

\vspace{5mm}
\noindent
T.H.~Baker$^{1}$ and  P.J.~Forrester$^{2}$

\noindent
$^{1}$Research Institute for Mathematical Sciences, 
Kyoto University, Kyoto 606, Japan; Present address
CMIS, 108 North Rd, Australian National University, Canberra 2601,
Australia

\noindent
$^{2}$Department of Mathematics and Statistics, University of Melbourne,  
Victoria 
3010, Australia; email: p.forrester@ms.unimelb.edu.au

\small
\begin{quote}
A bijection is given between fixed point free involutions of
$\{1,2,\dots,2N\}$ with maximum decreasing subsequence size
$2p$ and two classes of vicious (non-intersecting)
random walker configurations confined to the
half line lattice points
$l \ge 1$. In one class of walker configurations the
maximum displacement of the right most walker is $p$. Because the scaled
distribution of the  maximum decreasing subsequence size is known
to be in the soft edge GOE (random real symmetric matrices)
universality class, the same holds true for
the scaled distribution of the maximum displacement of the right most walker.
\end{quote}

Random permutations are fundamental combinatorial objects, which are
intimately related to other fundamental combinatorial objects
such as Young tableaux via the Robinson-Schensted-Knuth correspondence.
We recall that a Young tableau can be regarded as a numbered diagram of a
partition $\lambda_1 \ge \lambda_2 \ge \cdots \ge \lambda_p \ge 0$.
The diagram consists of squares drawn within a matrix array
with a square drawn in each row $(1 \le j \le p)$ and column
$k$ ($1 \le k \le \lambda_j$), while in each square is recorded a
number specified by some rule.
Recently random permutations, Young tableaux and their generalizations
have been shown to be at the core of certain statistical mechanical
models of growth processes \cite{Jo99a,Jo99b,PS99,GTW00}, 
vicious walker paths \cite{GOV98,Fo99,Ba00,KGV00} and
exclusion processes (the latter via mappings to certain growth
processes and vicious walker paths) 
amongst other topics. This has led to
progress in the study of these statistical mechanical models, by way
of the progress in the determination of fluctuation formulas for
quantities associated with random permutations \cite{BR00,BR99a,BR99b}.

As an example of the insight gained, we draw attention to the work
of Pr\"ahoffer and Spohn \cite{PS99}. These authors identify distinct
scaling forms for growth models in the 
Kardar-Parisi-Zhang (KPZ) universality class, that is growth models
described by the KPZ equation
$$
{\partial h \over \partial t} =
{\partial^2 h \over \partial x^2} + \Big (
{\partial h \over \partial x} \Big )^2 + \xi(t),
$$
where $\xi(t)$ is a noise term.
If the
growth profile is curved, the fluctuations are conjectured
to coincide with the distribution
of the largest eigenvalue in the scaled GUE (random Hermitian matrices),
while if the growth profile has zero curvature
the fluctuations are conjectured
to coincide with the distribution of the largest eigenvalue in the
scaled GOE (random real symmeric matrices). 
A matrix $X$ from the GOE $(\beta = 1)$ or GUE $(\beta = 2)$ is specified
by elements chosen with a joint distribution proportional to the
Gaussian $\exp(-\beta X^2/2)$. The largest eigenvalue occurs in the 
neighbourhood of $\lambda = \sqrt{2N}$
(which is referred to as the soft edge), and by making the scaling
$\lambda \mapsto  \sqrt{2N} + \lambda / \sqrt{2} N^{1/6}$ the
corresponding correlation functions have well defined limits \cite{Fo93a}.
Moreover, in the scaled $N \to \infty$ limit of both the GOE and GUE
the distribution of the largest eigenvalue is known
exactly  in terms of a certain Painlev\'e II
transcendent \cite{TW94a,TW96}. This identification was formulated from an exact
correspondence between a particular model of KPZ growth --- the
polynuclear growth model --- and increasing subsequences of random
permutations. The latter have been proved to have GUE soft edge
fluctuations \cite{BDJ98,Ok99,BOO99,Jo99b} 
in the absence of further constraints, but
GOE soft
edge fluctuations in the presence of the symmetry constraint restricting
the permutations to fixed point free involutions. The two cases
correspond to a curved and zero curvature interface respectively in the
corresponding polynuclear growth model.

In this work we will identify a statistical mechanical model
for which the profile displacement can be put into correspondence with
the maximum decreasing subsequence length of fixed point free involutions.
As this quantity has been shown rigorously to have GOE soft edge
fluctuations \cite{BR99b}, it follows that the profile of the
statistical mechanical model is in the GOE soft edge universality class.
We remark that the polynuclear growth model from a flat substrate is
also in correspondence with the maximum decreasing subsequence length
of fixed point free involutions \cite{PS99}, the mapping being quite
direct (unlike the present case). 
The model to be considered is the random turns model of vicious random
walkers. This model can be viewed either as a two-dimensional lattice
model of non-intersecting directed paths, or as a stochastic model of hard
core particles on a lattice in one-dimension. In the latter picture, at
discrete time intervals $t=1,2,\dots$ a particle which has a vacant site
as its left neighbour or its right neighbour (or both neighbours)
is selected at random and moved to the vacant
neighbouring site (if both sites are vacant, either is chosen with
equal probability). Plotting the trajectories of the particles on an
$l-t$ diagram ($l$ labelling the lattice sites)
gives the directed, non-intersecting paths picture
of the model.

Our interest is in two classes of configurations of this walker model.
The first is when there are exactly $p$ walkers, initially
equally spaced on neighbouring lattice sites $l=1,2,\dots,p$
and furthermore constrained to the region $l \ge 1$ (in the 
vicious walker vernacular, at the site $l=0$ there is a cliff at
which the walkers fall to their death \cite{Fo89e}). It is required that
after $2N$ steps the walkers return to their initial sites.

\begin{figure}
\epsfxsize=10cm
\centerline{\epsfbox{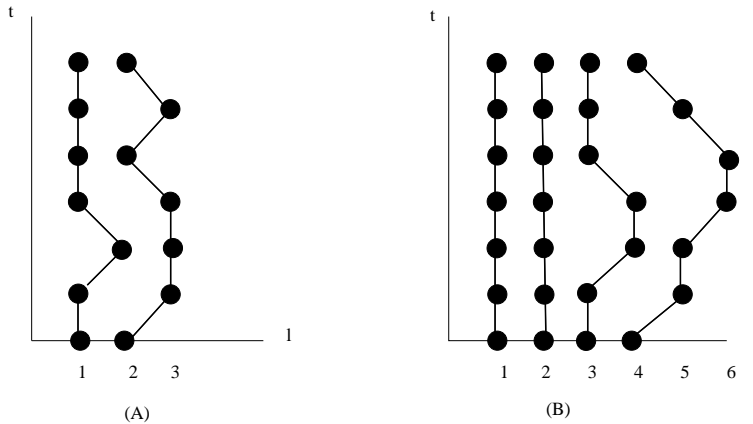}}
\caption{\label{f1} In diagram (A) there are $p=2$ walkers which move a
total of $2N=6$ steps. In diagram (B) there are $N^* = 4$ walkers which
move a total of $2N=6$ steps with the maximum displacement of the
rightmost walker given by $p=2$.
}
\end{figure}

In the second class of configurations the walkers again begin on
the neighbouring sites $l=1,2,\dots$, are confined to the region
$l \ge 1$, and return to their initial sites after $2N$ steps.
But rather than there being $p$ walkers there are now
$N^*$ walkers with $N^* \ge N$. The parameter $p$ enters by the
requirement that the right-most walker has a maximum displacement
of no more than $p$ lattice sites from its initial position. On the
other hand the value of $N^*$ is not a relevant parameter because
only a maximum  of $N$ consecutive walkers, counted from the
right-most walker, move from their initial sites. An example of
the first and second class of configurations is given in
Figure \ref{f1}.

Our first result is that both these classes of configurations are in
one-to-one correspondence with a certain subclass of fixed point
free involutions.

\begin{prop}\label{p1}
For both of the two classes of vicious walker 
configurations specified above there is a bijection with
fixed point free involutions of $\{1,2,\dots,2N\}$ (i.e.~permutations
consisting soley of two cycles) constrained so that the length
of the maximum decreasing subsequence is less than or equal to $2p$.
\end{prop}

We will give the details of the bijection for the second class of
configurations, and afterwards indicate the modification required
to establish the bijection for the first class of configurations.
Denote the walker initially at lattice site $j$ by $N^*+1-j$
$(j=1,2,\dots,N^*)$. A walker configuration can then be coded as a
sequence of integers from the alphabets $1,2,\dots,N^*$ and
$\bar{1}, \bar{2}, \dots, \bar{N}^*$, with the occurence of the
integer $j$ ($\bar{\jmath}$) at position $t$ denoting that
walker $j$ ($\bar{\jmath}$) moved one site to the right (left) at
time step $t$.

Of course not all words of length $2N$ from this alphabet give
rise to legal walker configurations. For a legal configuration,
at each time step we must have that
\begin{equation}\label{A}
n_1 \ge n_2 \ge \cdots \ge n_{N^*} \ge 0
\end{equation}
where
$$
n_j := \# j{\rm '}s - \# \bar{\jmath}{\rm '}s
$$
and after time step $2N$, each $n_j$ must equal zero. The requirement
(\ref{A}) can be represented diagramatically as the diagram of a
(conjugate) partition in which column $j$ is of length $n_j$
(see Figure \ref{f2}). A successive sequence of diagrams so generated
(starting and finishing with the empty diagram $\emptyset$)
uniquely specifies the walk thus demonstrating a bijection
between such diagrams and walk configurations.

In the theory of Young tableaux, the diagrams so generated are 
examples of oscillating tableaux. In general these tableaux map onto
certain two line arrays \cite{Su90}, which in the present case
represent involutions of $\{1,2,\dots,2N\}$ with no
fixed points. To construct the array, we number the box $i$ if
it is added to the diagram at step $i$. If instead a box is removed
at step $i$ (say from column $j$) this is to be done via the procedure
of reverse column insertion, which means if the particular box
ejected, $x_i$ say, was then inserted by the Schensted column
insertion procedure the original diagram would be restored
(see \cite{Fu97} for a description of the Schensted algorithm).
The fact that the removal occured at step $i$ is recorded by putting
the pair $(i,x_i)$ into a two line array with $i$ on top. Note
that since $x_i$ was bumped out at step $i$, it must have been
inserted in an earlier step, so $x_i < i$.
Furthermore all numbers in the array will be distinct, and at the
end of the procedure there will be $N$ pairs from the numbers
$\{1,2,\dots,2N\}$ with the top numbers ordered $i_1<i_2<\cdots < i_N$.
This procedure is illustrated in Figure \ref{f2}.
The pairs forming the array can be considered as the two cycles in
a fixed point free involution of $\{1,2,\dots,2N\}$.

\setlength{\unitlength}{4144sp}%
\begingroup\makeatletter\ifx\SetFigFont\undefined%
\gdef\SetFigFont#1#2#3#4#5{%
  \reset@font\fontsize{#1}{#2pt}%
  \fontfamily{#3}\fontseries{#4}\fontshape{#5}%
  \selectfont}%
\fi\endgroup%

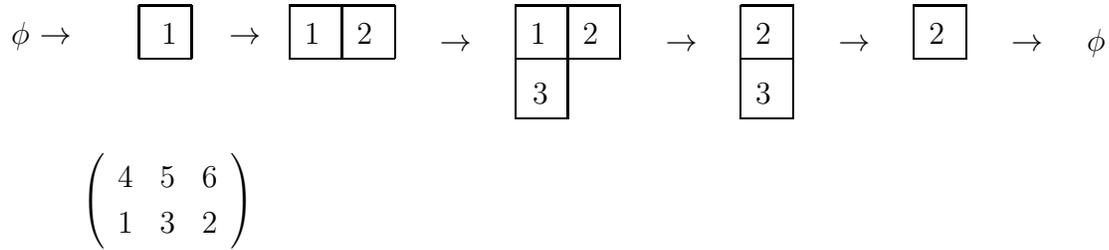
\begin{figure} 
\begin{picture}(6435,2295)(1036,-3571)
\thinlines
\put(2701,-2401){\line( 1, 0){315}}
\put(3016,-2401){\line( 0, 1){315}}
\put(3016,-2086){\line(-1, 0){315}}
\put(3016,-2401){\line( 1, 0){315}}
\put(3331,-2401){\line( 0, 1){315}}
\put(3331,-2086){\line(-1, 0){315}}
\put(1801,-2401){\line( 1, 0){315}}
\put(2116,-2401){\line( 0, 1){315}}
\put(2116,-2086){\line(-1, 0){315}}
\put(1801,-2086){\line( 0,-1){315}}
\put(2701,-2086){\line( 0,-1){315}}
\put(4051,-2401){\line( 1, 0){315}}
\put(4366,-2401){\line( 0, 1){315}}
\put(4366,-2086){\line(-1, 0){315}}
\put(4366,-2401){\line( 1, 0){315}}
\put(4681,-2401){\line( 0, 1){315}}
\put(4681,-2086){\line(-1, 0){315}}
\put(4051,-2086){\line( 0,-1){315}}
\put(4051,-2401){\line( 0,-1){360}}
\put(4051,-2761){\line( 1, 0){315}}
\put(4366,-2761){\line( 0, 1){360}}
\put(5401,-2401){\line( 1, 0){315}}
\put(5716,-2401){\line( 0, 1){315}}
\put(5716,-2086){\line(-1, 0){315}}
\put(5401,-2086){\line( 0,-1){675}}
\put(5401,-2761){\line( 1, 0){315}}
\put(5716,-2761){\line( 0, 1){360}}
\put(6436,-2401){\line( 1, 0){315}}
\put(6751,-2401){\line( 0, 1){315}}
\put(6751,-2086){\line(-1, 0){315}}
\put(6436,-2401){\line( 0, 1){315}}
\put(1396,-1411){\makebox(0,0)[lb]{\smash{\SetFigFont{12}{14.4}{\rmdefault}{\mddefault}{\updefault}
$1 2 1 \bar{2} \bar{1}  \bar{1}$ 
}}}
\put(1936,-2311){\makebox(0,0)[lb]{\smash{\SetFigFont{12}{14.4}{\rmdefault}{\mddefault}{\updefault}
1
}}}
\put(2341,-2311){\makebox(0,0)[lb]{\smash{\SetFigFont{12}{14.4}{\rmdefault}{\mddefault}{\updefault}
$\to$
}}}
\put(3601,-2356){\makebox(0,0)[lb]{\smash{\SetFigFont{12}{14.4}{\rmdefault}{\mddefault}{\updefault}
$\to$
}}}
\put(4951,-2356){\makebox(0,0)[lb]{\smash{\SetFigFont{12}{14.4}{\rmdefault}{\mddefault}{\updefault}
$\to$
}}}
\put(5986,-2356){\makebox(0,0)[lb]{\smash{\SetFigFont{12}{14.4}{\rmdefault}{\mddefault}{\updefault}
$\to$
}}}
\put(7021,-2356){\makebox(0,0)[lb]{\smash{\SetFigFont{12}{14.4}{\rmdefault}{\mddefault}{\updefault}
$\to$
}}}
\put(7471,-2356){\makebox(0,0)[lb]{\smash{\SetFigFont{12}{14.4}{\rmdefault}{\mddefault}{\updefault}
$\phi$
}}}
\put(1036,-2311){\makebox(0,0)[lb]{\smash{\SetFigFont{12}{14.4}{\rmdefault}{\mddefault}{\updefault}
$\phi$  $\to$
}}}
\put(1396,-3301){\makebox(0,0)[lb]{\smash{\SetFigFont{12}{14.4}{\rmdefault}{\mddefault}{\updefault}
$ \left ( \begin{array}{ccc} 4 & 5 & 6
\\ 1 & 3 & 2 \end{array} \right )$
}}}
\put(2791,-2311){\makebox(0,0)[lb]{\smash{\SetFigFont{12}{14.4}{\rmdefault}{\mddefault}{\updefault}
1
}}}
\put(3106,-2311){\makebox(0,0)[lb]{\smash{\SetFigFont{12}{14.4}{\rmdefault}{\mddefault}{\updefault}
2
}}}
\put(4141,-2311){\makebox(0,0)[lb]{\smash{\SetFigFont{12}{14.4}{\rmdefault}{\mddefault}{\updefault}
1
}}}
\put(4456,-2311){\makebox(0,0)[lb]{\smash{\SetFigFont{12}{14.4}{\rmdefault}{\mddefault}{\updefault}
2
}}}
\put(4096,-2671){\makebox(0,0)[lb]{\smash{\SetFigFont{12}{14.4}{\rmdefault}{\mddefault}{\updefault}
 3
}}}
\put(6526,-2311){\makebox(0,0)[lb]{\smash{\SetFigFont{12}{14.4}{\rmdefault}{\mddefault}{\updefault}
2
}}}
\put(5491,-2311){\makebox(0,0)[lb]{\smash{\SetFigFont{12}{14.4}{\rmdefault}{\mddefault}{\updefault}
2
}}}
\put(5491,-2671){\makebox(0,0)[lb]{\smash{\SetFigFont{12}{14.4}{\rmdefault}{\mddefault}{\updefault}
3
}}}
\put(1396,-3571){\makebox(0,0)[lb]{\smash{\SetFigFont{12}{14.4}{\rmdefault}{\mddefault}{\updefault}
}}}
\end{picture}
\caption{\label{f2} The word corresponding to the walker configuration (B) 
of Figure \ref{f1}, the sequence of oscillating tableaux
corresponding to the word, and the two line array constructed from the
oscillating tableaux.}
\end{figure}

The constraint that the rightmost walker have maximum displacement
of exactly $p$ lattice sites to the right of its starting
position means that the maximum length of the first column of
each tableau is less than or equal to
$p$ boxes. Because each tableaux has the number of
the boxes strictly increasing down each column and across each
row, it follows from the reverse column bumping procedure used to
form the corresponding two line array that the maximum
decreasing subsequence length in the bottom line is precisely
the maximum size of the first column (see the example of Figure
\ref{f2} ) which is less than or equal to $p$.
Hence the walker configurations are in one-to-one correspondence
with the two line arrays already noted subject to the additional
constraint that the maximum decreasing subsequence length in the
bottom line is less than or equal to $p$. In the
correspondence between the two line array and fixed point
free involutions, this constrains the fixed point free involutions
to have  maximum decreasing subsequence length less than or equal to
$2p$. To see this we note that the fixed point free involution can
be constructed by extending the top line of the two line array
to include all integers $1,2,\dots,2N$ in order and filling in the 
bottom line according to pairings implied by the original two line
array. We see that if
$$
x_{j_1} > x_{j_2} > \cdots > x_{j_q}
$$
is a particular decreasing subsequence of maximum length $q$ $(q \le p)$
in the bottom line of the original two line array, then the increasing
subsequence of length $2q$ formed from
$$
\{ j_1, j_2,\dots, j_q\} \cup \{ x_{j_1}, x_{j_2}, \dots, x_{j_q} \}
$$
in the top line of the new two line array gives a decreasing
subsequence of length $2q$ in the bottom line of the new two line
array. This construction worked in reverse shows that no decreasing
subsequence in the fixed point free involution can have length
greater than $2q$. 

The above procedure associating each walker configuration with a
two line array is reversible in that starting with a two line
array of the type specified a unique sequence of oscillating tableaux
and thus walker configuration can be constructed. Following
\cite{Su90} we work backwards
in the construction of the two line array 
from the sequence of oscillating tableaux. In going from the tableau at
step $i$ to that at step $i-1$  there are two distinct situations. One
is that $i$ does not appear in the top row of the two line array,
indicating that the tableau at step $i$ was not the result of removing
a box from the tableau at step $i-1$, but rather came from adding
a box labelled $i$ to the tableau at step $i-1$. Thus deleting the
box labelled $i$ from the tableau at step $i$ gives the tableau at
step $i-1$. On the other hand we may have that $i$ does appear in the
top row of the array, being part of the pair $(i,x_i)$. In this
case the tableau at step $i$ was obtained from the tableau at step
$i-1$ as a result of an inverse column bumping which ejected
$x_i$. Thus the tableau at step $i-1$ is constructed from the
tableau at step $i$ by Schensted column inserting $x_i$.
An example of this inverse procedure is given in
Figure \ref{f3}.
From the rules of the column insertion the maximum attained
length of the first column of the tableaux will equal the
length of the largest decreasing subsequence in the bottom line of the
two-line array and thus be less than or equal to $p$.

\setlength{\unitlength}{4144sp}%
\begingroup\makeatletter\ifx\SetFigFont\undefined%
\gdef\SetFigFont#1#2#3#4#5{%
  \reset@font\fontsize{#1}{#2pt}%
  \fontfamily{#3}\fontseries{#4}\fontshape{#5}%
  \selectfont}%
\fi\endgroup%
\begin{figure}
\begin{picture}(6030,1485)(1576,-2311)
\thinlines
\put(3601,-1411){\line( 0,-1){315}}
\put(3601,-1726){\line( 1, 0){315}}
\put(3916,-1726){\line( 0, 1){315}}
\put(3916,-1411){\line(-1, 0){315}}
\put(2611,-1411){\line( 0,-1){315}}
\put(2611,-1726){\line( 1, 0){315}}
\put(2926,-1726){\line( 0, 1){315}}
\put(2926,-1411){\line(-1, 0){315}}
\put(3916,-1411){\line( 0,-1){315}}
\put(3916,-1726){\line( 1, 0){315}}
\put(4231,-1726){\line( 0, 1){315}}
\put(4231,-1411){\line(-1, 0){315}}
\put(4861,-1411){\line( 0,-1){315}}
\put(4861,-1726){\line( 1, 0){315}}
\put(5176,-1726){\line( 0, 1){315}}
\put(5176,-1411){\line(-1, 0){315}}
\put(6661,-1411){\line( 0,-1){315}}
\put(6661,-1726){\line( 1, 0){315}}
\put(6976,-1726){\line( 0, 1){315}}
\put(6976,-1411){\line(-1, 0){315}}
\put(1576,-961){\makebox(0,0)[lb]{\smash{\SetFigFont{12}{14.4}{\rmdefault}{\mddefault}{\updefault}
$\left ( \begin{array}{ccc}3&4&6 \\ 1&2&5 \end{array} \right )$
}}}
\put(1621,-1636){\makebox(0,0)[lb]{\smash{\SetFigFont{12}{14.4}{\rmdefault}{\mddefault}{\updefault}
$\phi$
}}}
\put(2161,-1636){\makebox(0,0)[lb]{\smash{\SetFigFont{12}{14.4}{\rmdefault}{\mddefault}{\updefault}
$\leftarrow$
}}}
\put(3151,-1636){\makebox(0,0)[lb]{\smash{\SetFigFont{12}{14.4}{\rmdefault}{\mddefault}{\updefault}
$\leftarrow$
}}}
\put(4456,-1636){\makebox(0,0)[lb]{\smash{\SetFigFont{12}{14.4}{\rmdefault}{\mddefault}{\updefault}
$\leftarrow$
}}}
\put(5491,-1636){\makebox(0,0)[lb]{\smash{\SetFigFont{12}{14.4}{\rmdefault}{\mddefault}{\updefault}
$\leftarrow$
}}}
\put(5851,-1636){\makebox(0,0)[lb]{\smash{\SetFigFont{12}{14.4}{\rmdefault}{\mddefault}{\updefault}
$\phi$
}}}
\put(1621,-2311){\makebox(0,0)[lb]{\smash{\SetFigFont{12}{14.4}{\rmdefault}{\mddefault}{\updefault}
$12\bar{2}\bar{1}1\bar{1}$
}}}
\put(2746,-1636){\makebox(0,0)[lb]{\smash{\SetFigFont{12}{14.4}{\rmdefault}{\mddefault}{\updefault}
1
}}}
\put(3691,-1636){\makebox(0,0)[lb]{\smash{\SetFigFont{12}{14.4}{\rmdefault}{\mddefault}{\updefault}
1
}}}
\put(4051,-1636){\makebox(0,0)[lb]{\smash{\SetFigFont{12}{14.4}{\rmdefault}{\mddefault}{\updefault}
2
}}}
\put(4996,-1636){\makebox(0,0)[lb]{\smash{\SetFigFont{12}{14.4}{\rmdefault}{\mddefault}{\updefault}
2
}}}
\put(7606,-1636){\makebox(0,0)[lb]{\smash{\SetFigFont{12}{14.4}{\rmdefault}{\mddefault}{\updefault}
$\phi$
}}}
\put(7201,-1636){\makebox(0,0)[lb]{\smash{\SetFigFont{12}{14.4}{\rmdefault}{\mddefault}{\updefault}
$\leftarrow$
}}}
\put(6751,-1636){\makebox(0,0)[lb]{\smash{\SetFigFont{12}{14.4}{\rmdefault}{\mddefault}{\updefault}
5
}}}
\put(6256,-1636){\makebox(0,0)[lb]{\smash{\SetFigFont{12}{14.4}{\rmdefault}{\mddefault}{\updefault}
$\leftarrow$
}}}
\end{picture}
\caption{\label{f3} Correspondence between a two line array
corresponding to a fixed point free involution and a sequence
of oscillating tableaux, and the correspondence between the
oscillating tableaux and a word. The word is equivalent to
walker configuration (A) of Figure \ref{f1}, but translated at
least one lattice site to the right so that $N^* > N$, with
stationary walkers filling the intervening sites to the left
down to $l=1$. 
}
\end{figure}
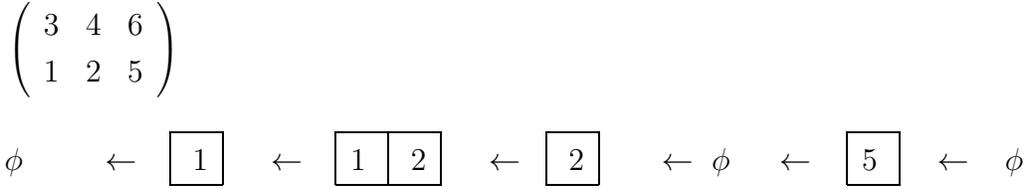

Hence for every sequence of oscillating tableaux, starting and
finishing with the empty tableau and having column length less
than or equal to $p$, there is a two line array equivalent
to a fixed point free involution
having maximum decreasing subsequence length less than or equal to
$2p$, and furthermore the correspondence can be established in the
reverse direction. Because there is a bijection between the
oscillating tableaux and random walker configurations, the result
of the Proposition \ref{p1} for the second class of walker configurations is
established.

Let us now turn our attention to the first class of configurations.
The walker configurations are again written as words, this time from
the alphabets $1,2,\dots,p$ and $\bar{1},\bar{2}, \dots, \bar{p}$. The
constraint (\ref{A}) (with $N^*$ replaced by $p$) is represented as
a diagram but now with row $j$ of length $n_j$ rather than column $j$
as previously, this feature being the essential difference between the
two cases. Note that the length of the first column now represents the
number of walkers displaced from their initial conditions.

As before, the boxes are numbered by $i$ if added at time step $i$,
and removed via the reverse column bumping procedure, with the
fact that the removal occured at step $i$ recorded by putting the
pair $(i,x_i)$ into a two line array with $i$ on top. An example is
given in Figure \ref{f4}.  The constraint that there
be less than or equal to $p$ walkers restricts the first column length
to be less than or equal to $p$. As already noted, the fact that the
reverse column bumping procedure is used to construct the two line array
from the tableau implies the former must therefore have decreasing
subsequence length no greater than $p$. Thus each walker configuration
can be mapped to a unique two line array of the same type as occured
in the corresponding mapping for the second class of configurations.
Futhermore, we have detailed how to associate such two line arrays
with a unique sequence of oscillating tableaux. From this sequence
of oscillating tableaux we can construct the word corresponding to
the walker configuration. Note that this differs from the construction
in the case of the second class of configurations because now it
is row $j$ which specifies the moves of walker $j$. The final result
is that there is a bijection between the first class of walker
configurations and two line arrays with top line ordered
$i_1<i_2< \cdots < i_N$, and maximum decreasing subsequence length no
greater than $p$. The latter being equivalent to fixed point free
involutions of $\{1,2,\dots,2N\}$ with maximum decreasing subsequence 
length no greater than $2p$, we see that Proposition \ref{p1} is now
established.

\setlength{\unitlength}{4144sp}%
\begingroup\makeatletter\ifx\SetFigFont\undefined%
\gdef\SetFigFont#1#2#3#4#5{%
  \reset@font\fontsize{#1}{#2pt}%
  \fontfamily{#3}\fontseries{#4}\fontshape{#5}%
  \selectfont}%
\fi\endgroup%
\begin{figure}

\setlength{\unitlength}{4144sp}%
\begin{picture}(5670,1587)(1666,-2986)
\thinlines
\put(3601,-1411){\line( 0,-1){315}}
\put(3601,-1726){\line( 1, 0){315}}
\put(3916,-1726){\line( 0, 1){315}}
\put(3916,-1411){\line(-1, 0){315}}
\put(2611,-1411){\line( 0,-1){315}}
\put(2611,-1726){\line( 1, 0){315}}
\put(2926,-1726){\line( 0, 1){315}}
\put(2926,-1411){\line(-1, 0){315}}
\put(3601,-1726){\line( 0,-1){315}}
\put(3601,-2041){\line( 1, 0){315}}
\put(3916,-2041){\line( 0, 1){315}}
\put(3916,-1726){\line(-1, 0){315}}
\put(4636,-1411){\line( 0,-1){315}}
\put(4636,-1726){\line( 1, 0){315}}
\put(4951,-1726){\line( 0, 1){315}}
\put(4951,-1411){\line(-1, 0){315}}
\put(6346,-1411){\line( 0,-1){315}}
\put(6346,-1726){\line( 1, 0){315}}
\put(6661,-1726){\line( 0, 1){315}}
\put(6661,-1411){\line(-1, 0){315}}
\put(2161,-1636){\makebox(0,0)[lb]{\smash{\SetFigFont{12}{14.4}{\rmdefault}{\mddefault}{\updefault}
$\to$
}}}
\put(3151,-1636){\makebox(0,0)[lb]{\smash{\SetFigFont{12}{14.4}{\rmdefault}{\mddefault}{\updefault}
$\to$
}}}
\put(2746,-1636){\makebox(0,0)[lb]{\smash{\SetFigFont{12}{14.4}{\rmdefault}{\mddefault}{\updefault}
1
}}}
\put(3691,-1636){\makebox(0,0)[lb]{\smash{\SetFigFont{12}{14.4}{\rmdefault}{\mddefault}{\updefault}
1
}}}
\put(3691,-1951){\makebox(0,0)[lb]{\smash{\SetFigFont{12}{14.4}{\rmdefault}{\mddefault}{\updefault}
2
}}}
\put(4231,-1636){\makebox(0,0)[lb]{\smash{\SetFigFont{12}{14.4}{\rmdefault}{\mddefault}{\updefault}
$\to$
}}}
\put(5221,-1636){\makebox(0,0)[lb]{\smash{\SetFigFont{12}{14.4}{\rmdefault}{\mddefault}{\updefault}
$\to$
}}}
\put(4726,-1636){\makebox(0,0)[lb]{\smash{\SetFigFont{12}{14.4}{\rmdefault}{\mddefault}{\updefault}
1
}}}
\put(1666,-2761){\makebox(0,0)[lb]{\smash{\SetFigFont{12}{14.4}{\rmdefault}{\mddefault}{\updefault}
$\left ( \begin{array}{ccc}3&4&6 \\
2 & 1 & 5 \end{array} \right )$
}}}
\put(5986,-1636){\makebox(0,0)[lb]{\smash{\SetFigFont{12}{14.4}{\rmdefault}{\mddefault}{\updefault}
$\to$
}}}
\put(6436,-1636){\makebox(0,0)[lb]{\smash{\SetFigFont{12}{14.4}{\rmdefault}{\mddefault}{\updefault}
5
}}}
\put(6931,-1636){\makebox(0,0)[lb]{\smash{\SetFigFont{12}{14.4}{\rmdefault}{\mddefault}{\updefault}
$\to$
}}}
\put(7336,-1636){\makebox(0,0)[lb]{\smash{\SetFigFont{12}{14.4}{\rmdefault}{\mddefault}{\updefault}
$\phi$
}}}
\put(1666,-2986){\makebox(0,0)[lb]{\smash{\SetFigFont{12}{14.4}{\rmdefault}{\mddefault}{\updefault}
}}}
\put(1846,-1636){\makebox(0,0)[lb]{\smash{\SetFigFont{12}{14.4}{\rmdefault}{\mddefault}{\updefault}
$\phi$
}}}
\put(5671,-1636){\makebox(0,0)[lb]{\smash{\SetFigFont{12}{14.4}{\rmdefault}{\mddefault}{\updefault}
$\phi$
}}}
\end{picture}
\caption{\label{f4} The oscillating tableaux for the configuration
(A) of Figure \ref{f1}, or equivalently the word of Figure \ref{f3},
constructed from the rules for a fixed $p$ of walkers, and 
the corresponding two line array.}
\end{figure}

The second class of configurations count the number of walker configurations
with a specific bound $p$ on the maximum displacement of the right-most
walker. From the $l-t$ diagrams of Figure \ref{f1} we see that this is
equivalent to counting the number of configurations which give rise to
a growth profile with a bound $p$ on its maximum spread, 
$L_N$ say. Our interest is in the distribution of $L_N$. Now,
with $L_t$ denoting the displacement of the right-most walker after
$t$ steps, the symmetry of the configurations under $t \to 2N - t$
means the statistical properties of $L_t$ are the same as those
of $L_{2N-t}$. At the centre of symmetry will be the
maximum displacement $L_N$, and the conjecture of
Pr\"ahoffer and Spohn predicts GUE fluctuations if the profile is
curved at this point, or GOE fluctuations if the profile has zero
curvature. Unfortunately the analytic form of the profile is not
known, so we cannot make use of this prediction presently.

In fact the nature of the fluctuations can be determined rigorously
by using the bijection of Proposition \ref{p1} between the second
class of configurations and fixed point free involutions. In
particular we can make use of the known distribution of the
maximum decreasing subsequence length for fixed point free involutions
of $\{1,2,\dots,2N\}$ to deduce the limiting distribution of $L_N$.
Regarding the former, let 
$L_N^{\rm inv}$ denote
the maximum decreasing subsequence length, and define the scaled
quantity
$$
\chi_N^{\rm inv} := {L_N^{\rm inv} - 2 \sqrt{2N} \over (\sqrt{2N})^{1/6} } =
{L_N^{\rm inv}/2 - 2 \sqrt{N/2} \over (\sqrt{N/2})^{1/6} }.
$$
Then it is proved in \cite{BR99b} that
$$
\lim_{N \to \infty} {\rm Pr} \Big ( \chi_N^{\rm inv} \le x \Big ) =
F_1(x)
$$
where $F_1(x)$ denotes the cumulative distribution of the largest
eigenvalue of matrices from the scaled GOE \cite{TW96,Fo99b}. The
following result is then an immediate consequence of
Proposition \ref{p1}.

\begin{prop}
Let $L_N$ denote the maximum displacement of the right-most
walker in the second class of random walker configurations specified
above, and set
$$
\chi_N := {L_N -  \sqrt{2N} \over {1 \over 2} (2N)^{1/6} }.
$$
Then
$$
\lim_{N \to \infty} {\rm Pr} \Big ( \chi_N \le x \Big ) =
F_1(x).
$$
\end{prop}

\noindent
Hence the walker profile at its maximum width
exhibits GOE fluctuations. The converse of the
prediction of Pr\"ahoffer and Spohn would then imply that the
walker profile has zero-curvature at this point.

As a final issue we consider the $p$-dimensional integral formula for the
number, $f_{Np}^{({\rm inv})}$ say, of fixed point free involutions of
$\{1,2,\dots,2N\}$ constrained so that the length of the maximum decreasing
subsequence is less than or equal to $2p$. With $USp(p)$ denoting the group
of $2p \times 2p$ unitary symplectic matrices (or equivalently the group
of $p \times p$ unitary matrices with real quaternion
elements ), it was shown by Rains \cite{Ra98} that
\begin{eqnarray}\label{rains}
f_{Np}^{({\rm inv})} & = & \Big \langle 
{\rm Tr}\, (S)^{2N} \Big \rangle_{S \in USp(p)} \nonumber \\
& = & {1 \over (2 \pi)^p p!}
\int_0^{ \pi} d \theta_1  \cdots
\int_0^{ \pi} d \theta_p 
\Big ( \sum_{j=1}^p 2 \cos \theta_j \Big )^{2N}
\prod_{j=1}^p |1 - z_j^2 |^2
\prod_{1 \le j < k \le p}
|1 - z_j z_k |^2 |z_j - z_k |^2,
\end{eqnarray}
where $z_j = e^{i \theta_j}$.

The formula (\ref{rains}) is in fact a special case of a counting formula
for a class of vicious walker paths. Thus consider
$p$ vicious walkers in the lock step model, confined to the lattice sites
$l \ge 1$, starting at positions
\begin{equation}\label{sn}
1 \le l_1^{(0)} < l_2^{(0)} < \cdots < l_p^{(0)},
\end{equation}
and arriving at positions
\begin{equation}\label{sn1}
1 \le l_1 < l_2 \cdots < l_p
\end{equation}
after $n$ steps. With $Z_n(l_1^{(0)}, \dots, l_p^{(0)};
 l_1, \dots, l_p)$ denoting the number of distinct walker configurations
of this prescription, we have the following result.

\begin{prop}
\begin{eqnarray}\label{2.9}
\lefteqn{Z_n(l_1^{(0)}, \dots, l_p^{(0)};
 l_1, \dots, l_p)} \nonumber \\
&& = {1 \over (2 \pi)^p}
\int_{-\pi}^\pi d \theta_1 \cdots \int_{-\pi}^\pi d \theta_p \,
\Big ( \sum_{j=1}^p 2 \cos \theta_j \Big )^n
\det \Big [ e^{i (l_j - l_k^{(0)})\theta_j } -
 e^{i (l_j + l_k^{(0)})\theta_j } \Big ]_{j,k=1,\dots,p}.
\end{eqnarray}
\end{prop}

Analogous to the proof of a similar counting formula in \cite{Fo91},
this can be verified by first noting from the definition of the
particular lock-step model that
$$
Z_n(l_1^{(0)}, \dots, l_p^{(0)};
 l_1, \dots, l_p) := Z_n(l_1, \dots, l_p)
$$
is the unique solution of the multidimensional difference equation
\begin{eqnarray}\label{2.10}
Z_{n+1}(l_1, \dots, l_p) & = &
Z_n(l_1 - 1, l_2, \dots, l_p) +
Z_n(l_1, l_2 - 1, \dots, l_p) \nonumber \\
&& + \cdots + Z_n(l_1,l_2,\dots, l_p - 1)  \nonumber \\
&&+ Z_n(l_1 + 1, l_2, \dots, l_p) +
Z_n(l_1, l_2 + 1, \dots, l_p) \nonumber \\
&& + \cdots + Z_n(l_1,l_2,\dots, l_p + 1)
\end{eqnarray}
subject to the non-intersection condition
\begin{equation}\label{2.11}
 Z_n(l_1, \dots, l_p) = 0 \quad {\rm if} \quad l_j = l_k \,\,
(j \ne k)
\end{equation}
the constraint $l_j \ge 1$ $(j=1,\dots,p)$ which requires
\begin{equation}\label{2.12}
 Z_n(l_1, \dots, l_p) = 0 \quad {\rm if} \quad l_1 = 0
\end{equation}
(here use has been make of the ordering (\ref{sn1})) and the initial
condition
\begin{equation}\label{2.13}
Z_n(l_1^{(0)}, \dots, l_p^{(0)};
 l_1, \dots, l_p) = \prod_{k=1}^p \delta_{l_k^{(0)}, l_k}
\end{equation}
where again use has been make of the orderings (\ref{sn}) and
(\ref{sn1}).

To verify that (\ref{2.9}) satisfies (\ref{2.10}) we note that (\ref{2.9})
gives
\begin{eqnarray}\label{2.14}
Z_{n+1}(l_1, \dots, l_p) & = & \Big ( {1 \over 2 \pi} \Big )^p
\sum_{\mu=1}^p \int_{-\pi}^\pi d \theta_1 \cdots
 \int_{-\pi}^\pi d \theta_p \,
(e^{i \theta_\mu} + e^{-i \theta_\mu}) \nonumber \\ &&
\times \Big ( \sum_{j=1}^p 2 \cos \theta_j \Big )^n
\det \Big [ e^{i (l_j - l_k^{(0)})\theta_j } -
 e^{i (l_j + l_k^{(0)})\theta_j } \Big ]_{j,k=1,\dots,p}
\end{eqnarray}
Using the fact that
\begin{eqnarray}
\lefteqn{
e^{\pm i \theta_\mu} \det
 [ e^{i (l_j - l_k^{(0)})\theta_j } -
 e^{i (l_j + l_k^{(0)})\theta_j } \Big ]_{j,k=1,\dots,p} } \nonumber \\
&& =
\det \left [  \begin{array}{c} e^{i (l_{j_1} - l_k^{(0)})\theta_j } -
 e^{i (l_{j_1} + l_k^{(0)})\theta_{j_1} } \nonumber \\
 e^{i (l_\mu \pm 1 - l_k^{(0)})\theta_\mu } -
 e^{i (l_\mu \pm 1 + l_k^{(0)})\theta_\mu }  \nonumber \\
 e^{i (l_{j_2} - l_k^{(0)})\theta_j } -
 e^{i (l_{j_2} + l_k^{(0)})\theta_j } 
\end{array}
\right 
]_{{j_1 = 1,\dots, \mu - 1 \atop j_2 = \mu+1, \dots,p} \atop k=1,\dots, p}
\end{eqnarray}
and recalling (\ref{2.9}) we can immediately identify the right hand side
of (\ref{2.13}) with the right hand side of (\ref{2.9}).

To verify (\ref{2.11}) we simply note that if $l_j = l_k$ for any
$j \ne k$ then two rows of the matrix in (\ref{2.9}) are the same so
the determinant vanishes. The condition (\ref{2.13}) is a property of 
(\ref{2.9}) since with $l_1 = 0$ the integrand
is odd in $\theta_1$ and thus the integral vanishes. Finally, to
verify the initial condition (\ref{2.13}) we make use of the
definition of a determinant
$$
\det [a_{jk}]_{j,k=1,\dots,p} = \sum_{P \in S_p} \varepsilon(P)
\prod_{j=1}^p a_{j P(j)},
$$
where $ \varepsilon(P)$ denotes the parity of the permutation $P$, to
expand the integrand in (\ref{2.9}) and integrate term by term.
Recalling each $l_j^{(0)}$ and $l_k$ is positive, this gives
\begin{equation}\label{2.18}
Z_0(l_1^{(0)}, \dots, l_p^{(0)};
 l_1, \dots, l_p) =
\sum_{P \in S_p} \varepsilon(P) \prod_{j=1}^p \delta_{l_j,
l_{P(j)}^{(0)}}.
\end{equation}
The ordering constraints (\ref{sn}) and (\ref{sn1}) imply that all terms
in (\ref{2.18}) except for the identity permutation must vanish, and so
(\ref{2.11}) is indeed satisfied.

Since the difference equation, the boundary conditions and the initial
conditions are all satisfied by (\ref{2.9}), we conclude that (\ref{2.9})
correctly represents $Z_n( l_1, \dots, l_p)$.

Although $Z_n$ is a positive integer, the integrand in 
(\ref{2.9}) is complex. An integral representation with a positive real
integrand can, in the case $l_j^{(0)} = l_j$ $(j=1,\dots,p)$, be
obtained by first noting
\begin{eqnarray}\label{A1}\lefteqn{
Z_{n_1+n_2}(l_1^{(0)}, \dots, l_p^{(0)};l_1^{(0)}, \dots, l_p^{(0)})}
\nonumber \\&&
= \sum_{1 \le l_1 < l_2 < \cdots < l_p }
Z_{n_1}(l_1^{(0)}, \dots, l_p^{(0)};l_1,\dots,l_p)
Z_{n_2}(l_1,\dots,l_p;l_1^{(0)}, \dots, l_p^{(0)}).
\end{eqnarray}
From (\ref{2.9}), and after simple manipulation of the determinant
therein, we see that
\begin{eqnarray}\label{A2}\lefteqn{
Z_{n_1}(l_1^{(0)}, \dots, l_p^{(0)};l_1,\dots,l_p)
Z_{n_2}(l_1,\dots,l_p;l_1^{(0)}, \dots, l_p^{(0)})} \nonumber \\&&
= {1 \over (2 \pi)^{2p}}
\int_{-\pi}^\pi d \theta_1 \cdots \int_{-\pi}^\pi d \theta_p
\int_{-\pi}^\pi d \phi_1 \cdots \int_{-\pi}^\pi d \phi_p
\Big (\sum_{j=1}^p 2 \cos \theta_j \Big )^{n_1}
\Big (\sum_{j=1}^p 2 \cos \phi_j \Big )^{n_2} \nonumber \\
&& \quad \times
\prod_{j=1}^p e^{i l_j(\theta_j - \phi_j)}
\det \Big [ e^{-i l_k^{(0)} \theta_j} - e^{i l_k^{(0)} \theta_j}
\Big ]_{j,k=1,\dots,p}
\det \Big [ e^{i l_k^{(0)} \phi_j} - e^{-i l_k^{(0)} \phi_j}
\Big ]_{j,k=1,\dots,p}
\end{eqnarray}
This is an even symmetric function of the $l_j$'s, which vanish for
$l_j = l_k$ and $l_j = 0$. Consequently the sum in (\ref{A1}) can be
replaced by
$$
{1 \over 2^p p! }  \sum_{l_1=-\infty}^\infty
 \sum_{l_2=-\infty}^\infty \cdots  \sum_{l_p=-\infty}^\infty
$$
Performing the sum in (\ref{A2}) using
$$
\sum_{l = - \infty}^\infty
e^{i l (\theta_j - \phi_j)} = 2 \pi \delta (\theta_j - \phi_j), \quad
|\theta_j - \phi_j| < 2 \pi
$$
and putting $n_1 + n_2 = n$ gives
\begin{eqnarray}\label{tut}
\lefteqn{Z_n(l_1^{(0)}, \dots, l_p^{(0)}; l_1^{(0)}, \dots, l_p^{(0)})}
\nonumber \\
&& = {1 \over (2 \pi)^p } {1 \over 2^p p!}
\int_{-\pi}^\pi d \theta_1 \cdots \int_{-\pi}^\pi d \theta_p \,
\Big (\sum_{j=1}^p 2 \cos \theta_j \Big )^{n}
\Big | \det \Big [ e^{i l_k^{(0)} \theta_j} -
e^{-i l_k^{(0)} \theta_j} \Big ]_{j,k=1,\dots,p} \Big |^2.
\end{eqnarray}

We can use (\ref{tut}) to rederive (\ref{rains}) since from the
definitions we have
\begin{equation}
f_{Np}^{\rm inv} = Z_{2N}(l_1^{(0)}, \dots, l_p^{(0)};
 l_1, \dots, l_p) \Big |_{l_j^{(0)} = l_j = j \, (j=1,\dots,p)} .
\end{equation}
Setting $l_j^{(0)} = j$ $(j=1,\dots,p)$ in (\ref{tut}), noting
that with $z_j = e^{i \theta_j}$ we have from the type C Vandermonde
formula \cite{Pr88}
\begin{equation}
\Big | \det [ z_j^k - z_j^{-k} ]_{j,k=1,\dots,p} \Big |^2 =
\prod_{j=1}^p |1 - z_j^2 |^2
\prod_{1 \le j < k \le p}
|1 - z_j z_k |^2 |z_j - z_k |^2,
\end{equation}
and making use of the fact that the integrand is even 
we see that the formula (\ref{rains}) indeed results.

\section*{Acknowledgements}
PJF acknowledges the financial support of the Australian Research Council,
and that of funds obtained by Prof.~K.Aomoto for his visit to Japan in
June 2000 which facilitated the present collaboration.


\end{document}